\documentclass[11pt,reqno]{amsart}
\usepackage{amsfonts}
\usepackage{amssymb}
\usepackage{graphicx}
\usepackage{float}
\usepackage{needspace}
\usepackage{amscd}
\usepackage[all]{xy}
\usepackage{amsmath,amsfonts,amssymb,amsthm,mathrsfs,pifont}
\usepackage[all]{xy}
\xyoption{curve}
\usepackage{fancybox}
\newtheorem{thm}{Theorem}[section]

\newtheorem{cor}[thm]{Corollary}

\newtheorem{lem}[thm]{Lemma}
\newtheorem{exm}[thm]{Example}

\newtheorem{prop}[thm]{Proposition}

\newtheorem{defn}[thm]{Definition}
\newtheorem*{remark*}{Remark}

\DeclareMathOperator{\Hod}{Hod}
\DeclareMathOperator{\Int}{Int}
\DeclareMathOperator{\gld}{gld}
\DeclareMathOperator{\q}{q}
\DeclareMathOperator{\rad}{rad}
\DeclareMathOperator{\soc}{soc}
\DeclareMathOperator{\End}{End}
\usepackage[hidelinks,bookmarks=false]{hyperref}
\hypersetup{pdftitle={Quasi-Hereditary Orderings of Nakayama Algebras},pdfauthor={Yuehui Zhang and Xiaoqiu Zhong}}
\title{Quasi-Hereditary Orderings of Nakayama Algebras}
\author{
Yuehui Zhang
\and
Xiaoqiu Zhong$^{*}$}
\address{School of Mathematical Sciences, Shanghai Jiao Tong University, Shanghai 200240, PR China}
\thanks{Supported by  National Natural Science Foundation of China under Grant No. 123B1026.}
\thanks{ * Corresponding author.}
\thanks{E-mail addresses: zyh@sjtu.edu.cn (Yuehui Zhang), 21-zxq@sjtu.edu.cn (Xiaoqiu Zhong)}
\thanks{Author's version of the article published in Bulletin of the Malaysian Mathematical Sciences Society (2026). The version of record is available at \url{https://doi.org/10.1007/s40840-026-02198-z}.}
\date{October 2026}

\begin{document}
\raggedbottom

\begin{abstract}
  Let $A$ be an algebra, and let $\mathcal{S}$ be the set of isomorphism classes of simple $A$-modules, with $|\mathcal{S}|=n$. In this paper, a total ordering of $\mathcal{S}$ is called a quasi-hereditary ordering, or a $q$-ordering, for $A$ if every Weyl module is Schurian and every indecomposable projective $A$-module is filtered by Weyl modules. The number of such total orderings is denoted by $\q(A)$. Determining whether a given total ordering of $\mathcal{S}$ is a $q$-ordering is a difficult problem. A result of Dlab and Ringel states that $A$ is hereditary if and only if every total ordering of $\mathcal{S}$ is a $q$-ordering; equivalently, $\q(A)=n!$. The $q$-ordering conjecture, proposed in 2000, asserts that every non-hereditary algebra $A$ satisfies $\q(A)\le\dfrac{2}{3}n!$. The main result of this paper is a necessary and sufficient criterion for a total ordering of the simple modules of a Nakayama algebra to be a $q$-ordering. More precisely, we show that a total ordering is a $q$-ordering if and only if, for every minimal relation $g$, the maximal element of $\Hod(g)$ does not belong to $\Int(g)$. As consequences, we characterize quasi-hereditary Nakayama algebras by the non-emptiness of the $X$-set, derive an iteration formula for $\q(A)$, and prove the $q$-ordering conjecture for Nakayama algebras. Finally, we give a non-monomial example showing that the conjecture does not hold for general finite-dimensional algebras.
\end{abstract}
\maketitle
\subsection*{Mathematics Subject Classification 2020}
16G20, 06A05, 16G10, 16D90
\subsection*{Keywords}
Quasi-hereditary algebras, Nakayama algebras, q-orderings, Weyl modules

\section{Introduction}

The theory of quasi-hereditary algebras has been extensively studied since its introduction by E. Cline, B. Parshall, and L. Scott in their seminal paper \cite{CPS} in 1988. The properties of a quasi-hereditary algebra $A$ are strongly dependent on a specific partial ordering of the set $\mathcal{S}$ of isomorphism classes of simple $A$-modules. This ordering ensures that the corresponding Weyl modules (or standard modules) are Schurian and that every indecomposable projective module admits a filtration by these modules. Such an ordering is called a \textit{quasi-hereditary ordering} for $A$. Consequently, a quasi-hereditary algebra is typically understood to come with a specified quasi-hereditary ordering.

It is a highly challenging problem to determine whether an algebra is quasi-hereditary. The first significant result was given by V. Dlab and C. M. Ringel \cite{DR89}, who showed that global dimension 2 implies quasi-heredity. The integral analogues of this result were obtained by S. K\"onig and A. Wiedemann \cite{K89} for quasi-hereditary orders. Thus, both hereditary algebras and hereditary orders are quasi-hereditary.

In 1991, C. M. Ringel \cite{ringer91} characterized quasi-hereditary algebras via a special tilting module called the characteristic module. Later, in 1996 and 1997, P. Zhang \cite{zp1, zp2} gave essential improvements to Ringel's characterization by weakening the characteristic module condition to the existence of a kind of faithful module. He also provided a criterion for the characteristic module by using a property of faithful and partial tilting modules.

Determining which orderings of $\mathcal{S}$ qualify as $q$-orderings for a given quasi-hereditary algebra remains an open challenge. Dlab and Ringel \cite{DR1992} found that two different partial orderings of $\mathcal{S}$ may yield the same Weyl modules. They defined an equivalence relation on quasi-hereditary partial orderings and proved that every such partial ordering is equivalent to a total ordering. Accordingly, throughout this paper a $q$-ordering means a total ordering of $\mathcal{S}$ for which $A$ is quasi-hereditary. We denote by $\mathscr{O}(A)$ the set of all such $q$-orderings for $A$.

This naturally leads to the question: How many $q$-orderings does a given algebra $A$ admit? Denote this number by $\q(A)$. Let $|\mathcal{S}| = n$ be the rank of the Grothendieck group of $A$. It is evident that $0 \leq \q(A) \leq n!$. A classical result, established by Dlab and Ringel, asserts that $A$ is hereditary if and only if $\q(A) = n!$---that is, if and only if every total ordering of $\mathcal{S}$ is a $q$-ordering \cite{DR89}.

Even two different total orderings may be equivalent. This equivalence relation relates to the concept of ``quasi-hereditary structures'' of an algebra $A$, denoted as \textbf{qh.st}($A$), which was defined and systematically studied by M. Flores, Y. Kimura, and B. Rognerud \cite{FKR}. Precisely, \textbf{qh.st}($A$) is the set of all equivalence classes of all quasi-hereditary partial orderings on $\mathcal{S}$. Clearly, \textbf{qh.st}($A$) may be equivalently defined as the set consisting of sets of Weyl modules of all $q$-orderings on $\mathcal{S}$. The notion of quasi-hereditary structures was originally mentioned by K. Coulembier \cite{C}. He proved that if an algebra $A$ admits a simple preserving duality, then $A$ has at most one quasi-hereditary structure.

Since a quasi-hereditary structure may contain several total orderings, classifying such structures and counting individual $q$-orderings are distinct problems. In this paper, we work with individual total orderings and their associated categories of \textit{good modules}.

The purpose of this paper is to provide a straightforward criterion for $q$-orderings of Nakayama algebras. To this end, we recall some known results about Nakayama algebras.

In 1990, M. Uematsu and K. Yamagata \cite{UY} demonstrated that Morita equivalence preserves quasi-heredity, meaning that if two algebras $A$ and $B$ are Morita equivalent, then $A$ is quasi-hereditary if and only if $B$ is quasi-hereditary. They also established that a non-hereditary Nakayama algebra $A$ is quasi-hereditary if and only if it has a simple module of projective dimension 2.

In 2000, the first author of the present paper and Y. Li \cite{ZL} proved that $\q(A) \leq \frac{2}{3} n!$ for tree-type quasi-hereditary algebras with a single generator. They conjectured that every non-hereditary algebra $A$ satisfies $\q(A) \leq \frac{2}{3}n!$; equivalently, $\q(A)>\frac{2}{3}n!$ if and only if $\q(A)=n!$, or, by the result of Dlab and Ringel, if and only if $A$ is hereditary. This statement is referred to as the ``$q$-ordering conjecture''. The final example of this paper shows that the conjecture does not hold for general non-monomial algebras; the positive result obtained here establishes it for Nakayama algebras. In 2008, the first author, along with L. Wu and C. Gao, proved the $q$-ordering conjecture when $A$ is an $\mathbb{A}$-type algebra with exactly two generators by deriving an explicit formula for $\q(A)$ in this context \cite{ZWG}.

In 2019, E. L. Green and S. Schroll \cite{GS} provided a necessary and sufficient condition for the existence of a $q$-ordering on $\mathcal{S}$ for a monomial algebra $A$ using quiver language. Their theory can be applied to the construction of concrete $q$-orderings by an inductive process on vertices of the underlying quiver of the involved algebra.

Among the homological approaches to quasi-heredity of Nakayama algebras, R. Marczinzik and E. Sen \cite{MS} introduced the following notion. A finite-dimensional algebra $A$ is called \textit{$S$-connected} if it has finite global dimension and the projective dimensions of its simple modules form an interval of nonnegative integers. They proved that a Nakayama algebra $A$ is $S$-connected if and only if $A$ is quasi-hereditary.

Y. Goto \cite{G} studied $q$-orderings obtained by permuting the simple modules in a given $q$-ordering. He gave a criterion for adjacent transpositions giving quasi-heredity, in terms of homological conditions of Weyl or Verma modules.

The homological criteria and the path-based approach of Green and Schroll address different aspects of quasi-heredity. Our aim is to test a prescribed total ordering directly from the minimal relations of a Nakayama algebra and to use this test to count all $q$-orderings. The criterion below avoids constructing the Weyl modules and their filtrations separately for each total ordering.

The main contribution of the present paper is Theorem \ref{quasi-hereditary order thm}, which gives a necessary and sufficient condition for a total ordering of the simple modules of a Nakayama algebra to be a $q$-ordering. Given a Nakayama algebra $A=KQ_n/I$, with $I=\langle g_1,\ldots,g_k\rangle$ minimally generated by paths, we show that a total ordering $\preceq$ is a $q$-ordering if and only if
\[
        \max \Hod(g_i)\notin \Int(g_i)
\]
for every generator $g_i$. Thus the problem of checking whether a given ordering is quasi-hereditary is reduced to a local combinatorial condition on the relations. The notions of hooks, denouements, hoods, and interiors are introduced precisely to state this condition in terms of these paths and their interiors.

Several consequences follow from this criterion. First, the quasi-heredity of a Nakayama algebra $A$ is determined by the subset $\mathscr{X}_A$, called the $X$-set of $A$: the algebra $A$ is quasi-hereditary if and only if $\mathscr{X}_A$ is nonempty (Theorem \ref{main thm}). The $X$-set can be read off directly from the interiors of the relations. The same criterion also gives an iteration formula for $\q(A)$ and the maximum of $\q(A)$ among non-hereditary Nakayama algebras with a fixed number of simple modules, thereby proving the $q$-ordering conjecture for Nakayama algebras (Theorem \ref{2/3}). The paper concludes with a separate non-monomial example showing that the conjecture fails for some finite-dimensional algebras outside the monomial setting.

Unless otherwise specified, all algebras discussed in this paper are assumed to be (connected) finite-dimensional over a fixed algebraically closed field $K$, with the unit element of an algebra $A$ denoted as $1_A$. For any set $X$, its cardinality is denoted by $|X|$.

\section{Preliminaries}
\subsection{Quasi-hereditary Algebra}
Let $A$ be an algebra, and let $\mathcal{S}$ be the set of isomorphism classes of simple $A$-modules. Fix a total ordering $\preceq$ of $\mathcal{S}$. For each $s\in\mathcal{S}$, choose a representative, again denoted by $s$, and let $P(s)$ be its projective cover. Define $\Delta(s)$ to be the largest factor module of $P(s)$ all of whose composition factors $t$ satisfy $t\preceq s$. The modules $\Delta(s)$, for $s\in\mathcal{S}$, are called Weyl modules.

Let ${\mathscr F}(\Delta)$ denote the class of all $A$-modules filtered by $\Delta$, meaning there exists a filtration
$$
0 = M_{t} \subset M_{t-1} \subset \cdots \subset M_{1} \subset M_{0} = M
$$
such that each factor module $M_{i-1}/M_{i}$ is isomorphic to some Weyl module, for all $1 \leq i \leq t$. Modules in ${\mathscr F}(\Delta)$ are referred to as good modules. For simplicity, ${\mathscr F}(\Delta)$ is denoted simply as ${\mathscr F}$.

\begin{defn}\label{2.1}
A total ordering $\preceq$ of $\mathcal{S}$ is called a \textit{quasi-hereditary ordering}, or a \textit{$q$-ordering}, for $A$ if, for every $s\in\mathcal{S}$,
\begin{enumerate}
    \item $\End_{A}(\Delta(s))\simeq K$;
    \item $P(s)\in{\mathscr F}$.
\end{enumerate}
The algebra $A$ is called quasi-hereditary if it admits a $q$-ordering.
\end{defn}

We denote by $\mathscr{O}(A)$ the set of all $q$-orderings for $A$, and write $\q(A)=|\mathscr{O}(A)|$.

A module $M$ satisfying $\End_{A}M\simeq K$ is called \textit{Schurian}. Thus a total ordering of $\mathcal{S}$ is a $q$-ordering for $A$ if and only if all the corresponding Weyl modules are Schurian and all indecomposable projective modules are good.

\subsection{Nakayama Algebras}

According to \cite{UY}, Proposition 2.3, Morita equivalence guarantees the equivalence of the involved module categories, hence preserves the property of being quasi-hereditary, as well as the number of quasi-hereditary orderings. Additionally, it is well-known that any finite-dimensional algebra is Morita equivalent to a quotient algebra of a path algebra defined by a finite quiver. Hence, we will concentrate our investigation on the quasi-heredity of Nakayama algebras by categorizing them into directed and cyclic types. Directed Nakayama algebras are those Morita equivalent to algebras defined by linearly oriented quivers, while cyclic Nakayama algebras are those Morita equivalent to algebras defined by cyclic quivers.

\subsubsection{Directed Nakayama Algebras}

\begin{figure}[H]
  \centering
  \includegraphics[scale=0.4]{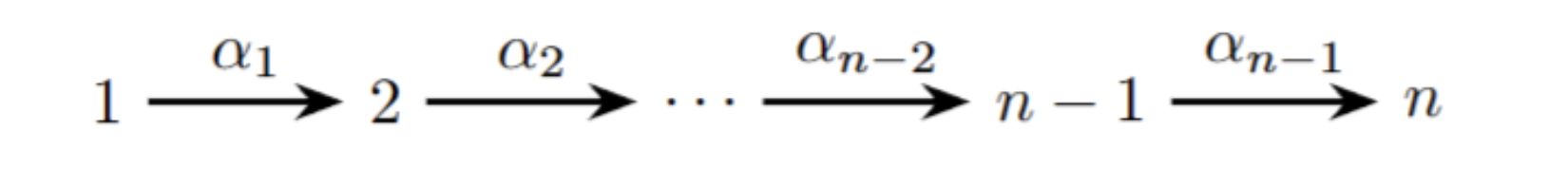}
  \caption{$A_n$}
  \label{fig:example1}
\end{figure}

For any $n\geq 1$, $A_n$ denotes the quiver shown in Figure \ref{fig:example1}.

The arrow from vertex $j$ to $j+1$ is denoted by $\alpha_j$, where $1\leq j\leq n-1$.

The oriented quiver $A_n$ is commonly referred to as the linearly oriented quiver of Dynkin type $A_n$.

An algebra $A$ is a directed Nakayama algebra of type $A_n$ if $A$ is a factor algebra $KA_n/I$ of the path algebra $KA_n$, where $I$ is an admissible two-sided ideal of $KA_n$.

\subsubsection{Cyclic Nakayama Algebra}

\begin{figure}[H]
  \centering
  \includegraphics{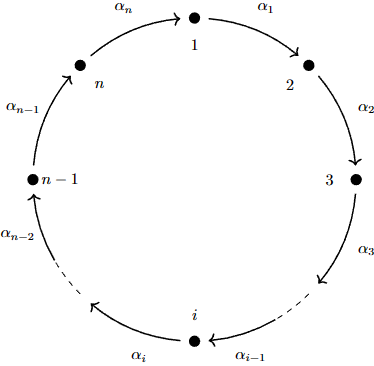}
  \caption{$\widetilde{A}_n$}
  \label{fig:example2}
\end{figure}

For $n\geq1$, let $\widetilde{A}_n$ denote the oriented cycle with $n$ vertices; for $n=1$, it consists of a single vertex with a loop $\alpha_1$. For $n\geq2$, it is depicted in Figure \ref{fig:example2}.

In the cyclic case, the arrows from vertex $j$ to $j+1$ and from vertex $n$ to $1$ are denoted by $\alpha_j$, where $1\leq j\leq n-1$, and $\alpha_n$ respectively.

An algebra $A$ is a cyclic Nakayama algebra if $A$ is a factor algebra $K\widetilde{A}_n/I$ of the path algebra $K\widetilde{A}_n$ where $I$ is an admissible two-sided ideal of $K\widetilde{A}_n$.

We use $Q_n$ to denote either $A_n$ or $\widetilde{A}_n$. For a Nakayama algebra $KQ_n/I$, the ideal $I$ is generated by paths. We fix a minimal set of path generators $\{g_1,\dots,g_k\}$ and write $I=\langle g_1,\dots,g_k\rangle$. These path generators are precisely the minimal zero relations, and below we simply call them relations. Each $g_i$ starts at a vertex $v_i$; after relabelling, we assume $v_1<v_2<\cdots<v_k$, with $k\leq n$. Such minimal sets of relations play a central role in the combinatorial and homological descriptions of cyclic Nakayama algebras; see \cite{RINGEL202195, Sen03062021}.

\subsection{Notations of Simple Modules, Hooks and Denouements of Paths}

Let $Q$ be any finite quiver, and let $A=KQ/I$ be a finite-dimensional algebra with $I$ an admissible ideal of $KQ$. For a vertex $x$, let $S_x$ denote the corresponding simple $A$-module, equivalently the simple top of the indecomposable projective module $P(x)$. We write $x$ for both the vertex and the isomorphism class of $S_x$. Thus the symbol $x$ may represent a number, a vertex, or an element of the set $\mathcal{S}$ of isomorphism classes of simple $A$-modules.

Now let $A=KQ_n/I$ be a Nakayama algebra. For a path generator $g$ in the fixed minimal generating set of $I$, the \textit{length} of $g$, denoted $\ell_g$, is the number of vertices it traverses (which equals the number of arrows plus one). The origin vertex of $g$ is referred to as the ``{\it hook}" of $g$, denoted by $h(g)$, while the end vertex of $g$ is termed the ``{\it denouement}" of $g$, denoted by $d(g)$. Thus $h(g)$ and $d(g)$ are respectively the source and the terminal vertex of the path $g$; the terminology is used to distinguish these two terms from the terms lying between them. The ordered multiset of simple $A$-modules encountered along $g$, including both endpoints and retaining repetitions, is called the ``{\it hood}" of $g$ and is denoted by $\Hod(g)$. The ordered multiset $\Int(g)$ is obtained from $\Hod(g)$ by deleting its first and last occurrences, corresponding to $h(g)$ and $d(g)$, respectively. In the cyclic case, the same simple module may occur more than once in $\Hod(g)$; hence, after the first and last occurrences have been deleted, the simple modules $h(g)$ or $d(g)$ may still belong to $\Int(g)$. Thus, $|\Hod(g)|=\ell_g$ and $|\Int(g)|=\ell_g-2$, where cardinalities count occurrences. Given any total ordering $\preceq$ of $\mathcal{S}$, denote by $\max \Hod(g)$ the maximum of $\Hod(g)$ with respect to $\preceq$. \footnote{For cyclic Nakayama algebras, $P(h(g))$ is a minimal projective module, $d(g)=\tau\soc P(h(g))$, and the denouements form the set $\mathcal S'$ used in \cite{Sen03062021}; see also \cite{RINGEL2013241}.}

\subsection{Definition of the $X$-set for Nakayama Algebras}
An arbitrarily given Nakayama algebra is Morita equivalent to $KQ_n/I$ for some $n$ with an admissible ideal $I$.

Since Morita equivalence preserves quasi-heredity, to analyze $q$-orderings of Nakayama algebras, it suffices to focus on Nakayama algebras of the form $KQ_n/I$. Each generator $g=\alpha_{d(g)-1}\alpha_{d(g)-2}\cdots\alpha_{h(g)}$ of $I$ is a path from its origin vertex $h(g)$ to its end vertex $d(g)$ of length $\ell_g\geq 3$, where $d(g)=h(g)+\ell_g-1$. In the cyclic case, vertex labels and arrow indices are read modulo $n$, using representatives in $\{1,\dots,n\}$, unless integer lifts are explicitly specified.

\begin{defn}
Let $A=KQ_n/I$ be a Nakayama algebra with $I=\langle g_1,\dots,g_k\rangle$, and put
\[
\mathcal H=\{h(g_i)\mid 1\leq i\leq k\},\qquad
\mathcal D=\{d(g_i)\mid 1\leq i\leq k\}.
\]
The set
\[
\mathscr{X}_A=\mathcal{S}\setminus\bigcup_{i=1}^{k}\Int(g_i)
\]
is called the \textit{$X$-set} of $A$. When no ambiguity is possible, we write $\mathscr{X}=\mathscr{X}_A$. We define three subsets of $\mathscr{X}$ separately by
\begin{align*}
\mathscr{X}^0&=\mathscr{X}\setminus(\mathcal H\cup\mathcal D),\\
\mathscr{X}^1&=\mathscr{X}\cap\bigl((\mathcal H\setminus\mathcal D)\cup(\mathcal D\setminus\mathcal H)\bigr),\\
\mathscr{X}^2&=\mathscr{X}\cap\mathcal H\cap\mathcal D.
\end{align*}
We also write $\mathscr{X}_A^j=\mathscr{X}^j$ for $j=0,1,2$.
\end{defn}

The sets $\mathscr{X}^0$, $\mathscr{X}^1$, and $\mathscr{X}^2$ are pairwise disjoint and form a partition of $\mathscr{X}$.

As we will see in Section 3 (Theorem \ref{main thm}), the quasi-heredity of a Nakayama algebra $A$ is completely determined by $\mathscr{X}_A$. The computation of $\q(A)$ is also highly dependent on $\mathscr{X}_A$.

\subsection{Notations for indecomposable modules over Nakayama Algebras}

For the basic Nakayama algebra $A=KQ_n/I$, each indecomposable module is uniserial and is determined up to isomorphism by its top and its $K$-dimension. We denote the indecomposable module with top $t$ and $K$-dimension $\ell$ by $(t;\ell)$, whenever it exists. In particular, $P(j)=(j;\ell)$ for $\ell=\dim_K P(j)$, and $\rad P(j)=(j+1;\ell-1)$ when $\ell>1$.

Condition (1) in Definition \ref{2.1} is automatically satisfied for directed Nakayama algebras. In the cyclic case, if $s$ is the largest simple module, then $\Delta(s)=P(s)$, so condition (1) requires $\dim_K P(s)\leq n$. The stronger necessary condition $\ell_{g_i}\leq n+1$ for every minimal relation in a quasi-hereditary Nakayama algebra will follow from Theorem \ref{quasi-hereditary order thm}.

Another easy fact about a directed Nakayama algebra is that it has a unique simple projective module $P(n)$ and $\Delta(n)=P(n)$ under every total ordering. Thus $P(n)$ always satisfies condition (2) in Definition \ref{2.1}.

Before delving into the general theory for Nakayama algebras, it is instructive to consider some concrete examples to illustrate the concept of $q$-orderings. The following simple cases demonstrate how the number of $q$-orderings, $\q(A)$, can be computed by directly verifying the conditions in Definition \ref{2.1}. These examples will also serve as a touchstone for the more powerful and general criterion we will develop later (Theorem \ref{quasi-hereditary order thm}).

\begin{exm}\label{A_3}
    Let $A=KA_3/\langle\alpha_2\alpha_1\rangle$ with the quiver shown in Figure \ref{fig:example3}.

\begin{figure}[H]
  \centering
  \includegraphics{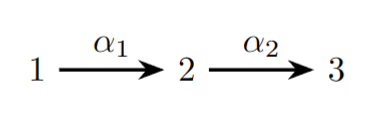}
  \caption{$A_3$}
  \label{fig:example3}
\end{figure}

The indecomposable projective modules are $P(1)=(1;2)$, $P(2)=(2;2)$, and $P(3)=(3;1)$. For the total ordering $1\succeq2\succeq3$, all composition factors of each $P(k)$ are at most its top, so $\Delta(k)=P(k)$ for $k=1,2,3$. Thus this is a $q$-ordering. For $1\succeq3\succeq2$, the maximal quotients defining the Weyl modules are $\Delta(1)=P(1)$, $\Delta(2)=2$, and $\Delta(3)=3$; the chain $P(2)\supset\rad P(2)=3\supset0$ is a $\Delta$-filtration, so this is also a $q$-ordering.

For $2\succeq1\succeq3$, the Weyl modules are $\Delta(1)=1$, $\Delta(2)=P(2)$, and $\Delta(3)=3$. A $\Delta$-filtration of $P(1)$ would require its radical $2$ to be filtered by these Weyl modules, which is impossible because the only Weyl module with top $2$ has dimension $2$. Hence this ordering is not a $q$-ordering. The ordering $2\succeq3\succeq1$ gives the same Weyl modules and also fails. For $3\succeq1\succeq2$, the Weyl modules are again $P(1),2,3$, while for $3\succeq2\succeq1$ they are all simple. Both orderings are therefore $q$-orderings, and $\q(A)=4$.
\end{exm}  

\begin{exm}\label{Tilde{A}3}
     Let $A=K\widetilde{A}_3/\langle\alpha_2\alpha_1\rangle$ with the quiver shown in Figure \ref{fig:example4}.

\begin{figure}[H]
  \centering
  \includegraphics{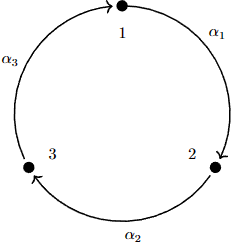}
  \caption{$\widetilde{A}_3$}
  \label{fig:example4}
\end{figure}

Similar to Example \ref{A_3}, one can easily get $\q(A)=4$. So the $q$-ordering conjecture holds in this case.
\end{exm} 

Marczinzik, Rubey, and Stump gave Kupisch-series criteria for simple modules of projective dimensions $1$ and $2$ \cite[Proposition 3.4]{MRS}. Sen's unpublished preprint \emph{Syzygy Filtrations of Cyclic Nakayama Algebras} (arXiv:1903.04645v3, 2022) studies syzygies in the cyclic case. Using hooks and denouements, we give criteria for these projective dimensions over Nakayama algebras, leading to a characterization of global dimension (denoted by $\gld$) greater than $2$.

\Needspace{12\baselineskip}
\begin{prop}\label{pd(s)=1}
    Let $A=KQ_n/I$ be a Nakayama algebra. In (1) and (2), let $s\in\mathcal{S}$ be a non-projective simple module. Then
\begin{enumerate}
    \item $pd(s)=1$ if and only if $s$ is not a hook.
    
    \item $pd(s)=2$ if and only if $s$ is a hook and its denouement belongs to $\mathscr{X}_A$.
    
    \item $\gld(A)\leq2$ if and only if $\mathcal D\subseteq\mathscr{X}_A$.
    
\end{enumerate}
\end{prop}

\begin{proof}
\begin{enumerate}
    \item $pd(s)=1$ if and only if $\rad P(s)=P(s+1)$ if and only if $s$ is not a hook.

    \item By (1), only a hook can have projective dimension $2$. Write $c_j=\dim_K P(j)$ and $e(j)=j+c_j$. In the cyclic case, use integer lifts of the vertex labels and extend $c_j$ periodically by $c_{j+n}=c_j$; thus $e(j+n)=e(j)+n$. Along each arrow, $c_{j+1}\geq c_j-1$, so $e(j+1)\geq e(j)$, with strict inequality precisely when $j$ is a hook.

    Suppose $s=h(g)$ and put $d=s+c_s$, an integer lift of $d(g)$. Since $\rad P(s)=(s+1;c_s-1)$, the kernel of its projective cover is
    \[
    \Omega^2(s)=(d;c_{s+1}-c_s+1).
    \]
    Here $\Omega^2(s)$ denotes the second syzygy of $s$, and it is nonzero because $s$ is a hook. Consequently,
    \begin{align*}
    pd(s)=2
    &\quad\Longleftrightarrow\quad c_d=c_{s+1}-c_s+1\\
    &\quad\Longleftrightarrow\quad e(d)=e(s+1)\\
    &\quad\Longleftrightarrow\quad\text{there is no hook among }s+1,\ldots,d-1.
    \end{align*}
    The last equivalence follows from the monotonicity of $e$ and its strict increase exactly at hooks.

    It remains to identify this condition with $d(g)\in\mathscr{X}_A$. If a hook $b$ satisfies $s<b<d$, then $e(b)\geq e(s+1)>e(s)=d$. Hence $d$ lies in the interior of the relation starting at $b$. Conversely, if $d$ lies in the interior of a relation starting at a hook $b$, then $b<d<e(b)$. Monotonicity excludes $b\leq s$, because that would give $e(b)\leq e(s)=d$. Thus $s<b<d$. In the cyclic case, membership of the vertex $d(g)$ in an interior is represented by such an occurrence after shifting the relation by a multiple of $n$. This proves (2), including relations that traverse a vertex more than once.

    \item The global dimension is the maximum of the projective dimensions of the simple modules, allowing the value $\infty$. Projective simples have projective dimension $0$, and non-hook, non-projective simples have projective dimension $1$ by (1). By (2), all hook simples have projective dimension $2$ if and only if all their denouements belong to $\mathscr{X}_A$. Thus $\gld(A)\leq2$ if and only if $\mathcal D\subseteq\mathscr{X}_A$.

\end{enumerate}

  \end{proof}

Proposition \ref{pd(s)=1} expresses the Kupisch-series criteria of Marczinzik, Rubey, and Stump \cite[Proposition 3.4]{MRS} in terms of relations. Write $c_s=\dim_KP(s)$, using our labels $1,\ldots,n$ in place of their labels $0,\ldots,n-1$, with the same order of composition factors. For a non-projective simple module, their condition $c_{s+1}+1=c_s$ for projective dimension $1$ says exactly that $s$ is not a hook. For a hook $s=h(g)$, their condition for projective dimension $2$ is
\[
c_{s+1}+1=c_s+c_{s+c_s}.
\]
Here $s+c_s$ is the denouement of $g$, with cyclic indices read modulo $n$. The proof above identifies this equality with $d(g)\in\mathscr X_A$. Applying these tests to all non-projective simples gives the criterion $\mathcal D\subseteq\mathscr X_A$ for global dimension at most $2$, in agreement with \cite[Theorem 4.2]{MRS}.

According to \cite{UY}, a non-hereditary Nakayama algebra $A$ is quasi-hereditary if and only if it has a simple module of projective dimension $2$. Now, one can use Proposition \ref{pd(s)=1} to quickly identify all simple modules of projective dimension $2$.

\Needspace{12\baselineskip}
\begin{exm}\label{A_5}
    Let $A=KA_5/\langle\alpha_2\alpha_1,\alpha_3\alpha_2\rangle$ with the quiver shown in Figure \ref{fig:C}.
    \begin{figure}[H]
        \centering
        \includegraphics{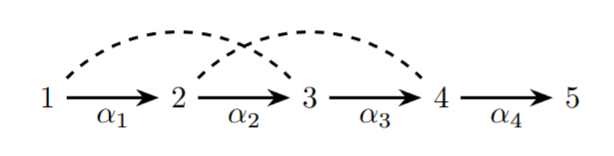}
        \caption{$A_5$}
        \label{fig:C}
    \end{figure}
    The Kupisch series is $(c_1,c_2,c_3,c_4,c_5)=(2,2,3,2,1)$, and
    \[
    \mathcal H=\{1,2\},\qquad \mathcal D=\{3,4\},\qquad
    \mathscr X_A=\{1,4,5\}.
    \]
    By Proposition \ref{pd(s)=1}, $pd(3)=pd(4)=1$ because $3$ and $4$ are not hooks, and $pd(2)=2$ because $2=h(\alpha_3\alpha_2)$ and $d(\alpha_3\alpha_2)=4\in\mathscr X_A$. In the Kupisch-series description, $c_{s+1}+1=c_s$ holds at $s=3,4$, while the test at the two hooks gives
    \[
    c_3+1=4=c_2+c_4,\qquad c_2+1=3\neq c_1+c_3=5.
    \]
    These agree with $4\in\mathscr X_A$ and $3\notin\mathscr X_A$, respectively. Since $\rad P(1)=2$, we have $pd(1)=1+pd(2)=3$. Since $5$ is projective, the projective dimensions are $(3,2,1,1,0)$ and $\gld(A)=3$. This same Kupisch series occurs in \cite[Figure 2]{MRS}; shifting its vertex labels by $1$ gives the same calculation. In the Dyck-path description of that paper, the tests for projective dimensions $1$ and $2$ correspond to double falls and rectangles, respectively.
\end{exm}

\section{Quasi-hereditary  Nakayama Algebras}
In this section, we establish a criterion for $q$-orderings over Nakayama algebras through the use of hoods. Subsequently, we derive a concise characterization for the quasi-heredity of Nakayama algebras, where the non-emptiness of the $X$-set serves as a necessary and sufficient condition.

\subsection{Quasi-hereditary orderings of Nakayama Algebras}

 Let $g$ be a path generator of $I$. If the cardinality  $|\Hod (g)|$ of the hood  $\Hod (g)$  is sufficiently small, then the indecomposable projective module $P(h(g))$ is Schurian, and hence the Weyl module $\Delta(h(g))$ is also Schurian, as stated in the following lemma.
 
\begin{lem}\label{Schurian lemma}
Let $A=KQ_n/I$ be a Nakayama algebra with $I=\langle g_1,\dots,g_k\rangle$. If $|\Hod (g_i)|\leq n+1$, $i=1,\dots,k$, then every Weyl module $\Delta(h(g_i))$ is Schurian under any total ordering of $\mathcal{S}$.
\end{lem}

\begin{proof}
For each minimal generator $g_i$, the indecomposable projective module $P(h(g_i))$ has top $h(g_i)$ and $K$-dimension $\ell_{g_i}-1$. The hypothesis $|\Hod(g_i)|=\ell_{g_i}\leq n+1$ gives $\ell_{g_i}-1\leq n$, so its consecutive composition factors are distinct and $P(h(g_i))$ is multiplicity free. Since $\Delta(h(g_i))$ is a quotient of $P(h(g_i))$, it is a multiplicity-free uniserial module. Any endomorphism of a multiplicity-free uniserial module is determined by its action on the top simple module, and a non-zero endomorphism cannot have image in the radical because no simple composition factor is repeated. Hence $\End_A(\Delta(h(g_i)))\simeq K$, and $\Delta(h(g_i))$ is Schurian.
\end{proof}

Hooks and denouements play key roles in the following result.

\begin{lem}\label{Ladder lemma}
Let $A=KQ_n/I$ be a Nakayama algebra, $1\leq k\leq n$. Given any total ordering $\preceq$ of $\mathcal{S}$ with good module category $\mathscr{F}$, the following statements hold.

(1) If  $\rad P(k)\in\mathscr{F}$, then $P(k)\in\mathscr{F}$.

(2) If $k$ is not a hook, there is an arrow from $k$ to $k+1$, and $P(k+1)\in \mathscr{F}$, then $P(k)\in \mathscr{F}$.

(3) Let $g$ be a generator of $I$ with hook $h(g)=k$ and denouement $d(g)$. If $d(g)=\max \Hod (g)$ and $P(k)$ is Schurian, then $P(k)\in \mathscr{F}$.
\end{lem}

\begin{proof}

(1) Write $P(k)=(k;\ell)$ and $\Delta(k)=(k;\ell')$, where $1\leq\ell'\leq\ell$. If $\ell'=\ell$, then $P(k)=\Delta(k)\in\mathscr F$. Otherwise, put $M=(k+\ell';\ell-\ell')$, so that $P(k)/M=\Delta(k)$.

By maximality of $\Delta(k)$, $k\preceq k+\ell'$ and $k\neq k+\ell'$, whereas $i\preceq k$ for every $i\in\{k,k+1,\dots,k+\ell'-1\}$. Thus no Weyl module whose top occurs before $M$ in $\rad P(k)$ can contain the simple module $k+\ell'$. Since $\rad P(k)$ is uniserial, any $\Delta$-filtration of it must therefore have $M$ as a filtration term. This includes $\ell'=1$, when $M=\rad P(k)$. Hence $M\in\mathscr F$, and adjoining the quotient $\Delta(k)$ gives a $\Delta$-filtration of $P(k)$.

(2) If $k$ is not a hook, then $\rad P(k)=P(k+1)$. So $\rad P(k)\in \mathscr{F}$. By (1), $P(k)\in \mathscr{F}$.

(3) Put $m=\ell_g-1=\dim_K P(k)$. In the cyclic case, the Schurian hypothesis gives $m\leq n$, since otherwise the top of $P(k)$ would repeat and yield a nonzero endomorphism with image in its radical. Thus, in either case, the composition factors $k,k+1,\ldots,k+m-1$ of $P(k)$ are distinct, with cyclic labels read modulo $n$.

Put $a_0=0$. Whenever $a_r<m$, define
\[
a_{r+1}=\min\bigl(\{b\mid a_r<b<m,\ k+a_r\preceq k+b\}\cup\{m\}\bigr).
\]
Thus we cut before the first composition factor strictly larger than the current top, or at the end of $P(k)$ if there is none. This gives $0=a_0<a_1<\cdots<a_t=m$.

For each $r<t$, all vertices $k+a_r,\ldots,k+a_{r+1}-1$ are at most $k+a_r$. Since $g$ is a minimal zero relation, every proper subpath of $g$ is nonzero in $A$, so $P(k+a_r)$ has an initial segment with these composition factors. Any larger quotient contains $k+a_{r+1}$, which is strictly larger than $k+a_r$ by the choice of the cut when $a_{r+1}<m$, and by the maximality of $k+m=d(g)$ when $a_{r+1}=m$ and $k+a_r\neq d(g)$. The only possible equality $k+a_r=d(g)$ occurs when $a_r=0$ and $h(g)=d(g)$; then $P(k)$ itself has dimension $m$, so no larger quotient exists. By the definition of the Weyl module,
\[
\Delta(k+a_r)=(k+a_r;a_{r+1}-a_r).
\]
Consequently,
\[
P(k)=\rad^{a_0}P(k)\supset\rad^{a_1}P(k)\supset\cdots\supset\rad^{a_t}P(k)=0
\]
has successive factors $(k+a_r;a_{r+1}-a_r)=\Delta(k+a_r)$, so it is a $\Delta$-filtration of $P(k)$. Here $\rad^aP(k)$ denotes the $a$th radical of $P(k)$, with $\rad^0P(k)=P(k)$.
\end{proof}

Now we can prove the following criterion for $q$-orderings. 

\begin{thm}\label{quasi-hereditary order thm}
Let $A=KQ_n/I$ be a Nakayama algebra with $I=\langle g_1,\dots,g_k\rangle$. Let $\preceq$ be a total ordering of $\mathcal{S}$. Then $\preceq$ is a quasi-hereditary ordering if and only if $\max \Hod (g_i)\notin \Int(g_i)$ for every $i=1,\dots,k$.
\end{thm}

\begin{proof}
Necessity. Assume that $\preceq$ is a $q$-ordering. We show that no hood has its maximum in its interior. Suppose, on the contrary, that there exists a relation $g$ such that the maximum $o=\max\Hod(g)$ belongs to $\Int(g)$. There are two cases.

(1) If $o=h(g)$ or $o=d(g)$ as a simple module, then the same simple module still belongs to $\Int(g)$ after the first and last terms of $\Hod(g)$ have been deleted. If $o=h(g)$, then $\Delta(o)=P(o)$, and its top and an interior occurrence give two occurrences of $o$. If $o=d(g)\neq h(g)$, the suffix from an interior occurrence of $o$ to the endpoint is a proper subpath of $g$, hence is nonzero by minimality. All its vertices are at most $o$, so $\Delta(o)$ again contains $o$ at least twice. A repeated top in a uniserial Nakayama module gives a nonzero endomorphism with image in the radical. Hence $\dim_K\End_{A}(\Delta(o))\ge 2$, so condition (1) in Definition \ref{2.1} fails, contradicting that $\preceq$ is a $q$-ordering.

(2) Suppose that $o\neq h(g),d(g)$ as a simple module. Put $m=\dim_K P(h(g))=\ell_g-1$, and let $t$ be the position of the first occurrence of $o$ among the composition factors of $P(h(g))$, numbered from $0$ to $m-1$. Thus $1\leq t\leq m-1$. Every composition factor preceding this occurrence is strictly smaller than $o$. In any $\Delta$-filtration of the uniserial module $P(h(g))$, the factor containing this occurrence must therefore begin at position $t$ and be isomorphic to $\Delta(o)$, since a Weyl module with smaller top cannot contain $o$.

The suffix of $g$ from position $t$ to its endpoint is a proper subpath of $g$, so it is nonzero by minimality of $g$. All its vertices are at most $o$. Hence $\dim_K\Delta(o)\geq m-t+1$. However, only $m-t$ composition factors remain in $P(h(g))$ from position $t$ onwards. This contradiction shows that $P(h(g))$ cannot admit a $\Delta$-filtration, contrary to condition (2) in Definition \ref{2.1}.

Sufficiency. Assume that $\preceq$ satisfies $\max \Hod(g_i)\notin \Int(g_i)$ for every $i=1,\dots,k$. Then $|\Hod(g_i)|\leq n+1$ and $\max \Hod(g_i)$ is either $h(g_i)$ or $d(g_i)$ for every $i$.

By Lemma \ref{Schurian lemma}, $\Delta(s)$ is Schurian whenever $s$ is a hook. Let $s$ be a non-hook vertex. Since $\Delta(s)$ is uniserial, it suffices to rule out a second occurrence of its top $s$. Such a repetition is possible only in the cyclic case. Its first $n$ composition factors would then contain every simple module, so $s$ would be the maximum of $\mathcal S$ and $\dim_K P(s)>n$.

Finite-dimensionality of the cyclic algebra guarantees at least one hook. Let $b$ be the first hook encountered after $s$ along the cycle, at distance $t$, where $1\leq t<n$. Repeatedly using $\rad P(j)=P(j+1)$ at the intervening non-hook vertices gives $\dim_K P(s)=t+\dim_K P(b)$. The minimal relation $g$ starting at $b$ has $\dim_K P(b)$ arrows. Measured from $s$, this occurrence of $g$ starts at position $t$ and ends at position $t+\dim_K P(b)$. Since $t<n<t+\dim_K P(b)$, the vertex $s$ at position $n$ belongs to $\Int(g)$. But $s=\max\Hod(g)$, contradicting the assumed hood condition. Thus $\Delta(s)$ is Schurian, and condition (1) of Definition \ref{2.1} holds.

It remains to verify condition (2) of Definition \ref{2.1}. Fix a relation $g_i$. If $h(g_i)=\max\Hod(g_i)$, then every composition factor of $P(h(g_i))$ is at most $h(g_i)$, and hence $\Delta(h(g_i))=P(h(g_i))$. Thus $P(h(g_i))\in\mathscr{F}$. Otherwise, $d(g_i)$ is the maximum. The projective $P(h(g_i))$ is multiplicity free, hence Schurian, because $\dim_K P(h(g_i))=\ell_{g_i}-1\leq n$. Lemma \ref{Ladder lemma}(3) then gives $P(h(g_i))\in\mathscr{F}$. Consequently, the projective module associated with every hook is good.

We now consider the projective modules associated with non-hook vertices. Between two consecutive hooks, start with the good projective module at the latter hook and move backwards along the quiver. At each non-hook vertex $j$, we have $\rad P(j)=P(j+1)$; hence Lemma \ref{Ladder lemma}(2) shows successively that $P(j)$ is good. In the cyclic case, the intervals between consecutive hooks cover all non-hook vertices. In the directed case, the same argument applies between consecutive hooks and before the first hook; after the last hook, it starts with the simple projective module $P(n)=\Delta(n)$ and again proceeds backwards. If there are no hooks, finite-dimensionality forces $Q_n=A_n$, and the same backwards induction from $P(n)$ covers the whole quiver. Therefore every indecomposable projective module belongs to $\mathscr{F}$.

This finishes the proof.
\end{proof} 

Theorem \ref{quasi-hereditary order thm} gives a criterion to determine which total orderings of $\mathcal{S}$ are $q$-orderings for Nakayama algebras. Specifically, it suffices to check the restriction of the ordering to each hood of $A$: the maximal simple module of each hood must not belong to its interior. Consequently, the verification of $q$-orderings is reduced to examining the corresponding total orderings for one-generator Nakayama algebras.

 We apply Theorem \ref{quasi-hereditary order thm} to check Example \ref{A_3}.

$A$ has only one generator $g=\alpha_2\alpha_1$ and $\Hod(g)=\Hod(\alpha_2\alpha_1)=\{1,2,3\}$, $\Int(\alpha_2\alpha_1)=\{2\}$. By Theorem \ref{quasi-hereditary order thm}, a total ordering $\preceq$ is quasi-hereditary if and only if $2$ is not the maximum of $\{1,2,3\}$ under $\preceq$. Therefore, $2\succeq 1\succeq 3$ and $2\succeq 3\succeq 1$ are not quasi-hereditary, while the other four total orderings are quasi-hereditary.

Using Theorem \ref{quasi-hereditary order thm} and the $X$-set, we can now prove the following criterion to determine whether a Nakayama algebra is quasi-hereditary. 

\begin{thm}\label{main thm}
    A Nakayama algebra $A$ is quasi-hereditary if and only if $\mathscr{X}_A\neq\emptyset$ .
     \begin{proof}

    Necessity. Suppose $\mathscr{X}_A=\emptyset$ and choose $\preceq\in\mathscr{O}(A)$. Let $s$ be the maximum of $\mathcal{S}$ under $\preceq$. Since $s\notin\mathscr{X}_A$, it belongs to $\Int(g_i)$ for some relation $g_i$. As the maximum of $\mathcal{S}$, it is also the maximum of $\Hod(g_i)$, contradicting Theorem \ref{quasi-hereditary order thm}.

Sufficiency. Let $A=KQ_n/I$ with $I=\langle g_1,\dots,g_k\rangle$. Since $\mathscr{X}_A\neq\emptyset$, $\ell_{g_i}\leq n+1$, 
$i=1,\dots,k$. Thus $d(g_i)\notin \Int(g_i)$, $i=1,\dots,k$. Let $x\in \mathscr{X}_A$. Consider the following total ordering $\preceq$:
\begin{align}\label{Reverse order}
    x+1\preceq x+2\preceq\dots\preceq n\preceq 1\preceq\dots\preceq x-1\preceq x
\end{align}
Since $x\notin\Int(g_i)$, the segment from the first interior vertex of $g_i$ to $d(g_i)$ does not cross the cut immediately after $x$ in the cyclic listing underlying ordering (\ref{Reverse order}). Thus every interior vertex strictly precedes $d(g_i)$, and $\max \Int(g_i)\preceq d(g_i)$ for $i=1,\dots,k$. Namely, $\max \Hod(g_i)\notin \Int(g_i)$ for every $i=1,\dots,k$. By Theorem \ref{quasi-hereditary order thm}, the total ordering $\preceq$ defined by (\ref{Reverse order}) is a $q$-ordering.
    \end{proof}
\end{thm}

Theorem \ref{main thm} gives a relation-based characterization of quasi-heredity for Nakayama algebras: once the relations are known, one only has to check whether some simple module does not occur in the interior of any hood.

We apply Theorem \ref{main thm} to characterize the quasi-heredity of Nakayama algebras associated with the quiver $\widetilde{A}_3$. The key cases are summarized below:

\begin{description}
    \item[Case 1] $K\widetilde{A}_3/\langle\alpha_3\alpha_2\alpha_1\rangle$ \\
          $\mathscr{X}=\{1\}$ (i.e., $\mathscr{X}\neq\emptyset$), so $A$ is quasi-hereditary by Theorem \ref{main thm}.
    \item[Case 2] $K\widetilde{A}_3/\langle\alpha_1\alpha_3\alpha_2\alpha_1\rangle$ \\
          $\mathscr{X}=\emptyset$, so $A$ is not quasi-hereditary by Theorem \ref{main thm}. Here the single relation has length 5.
    \item[Case 3] $K\widetilde{A}_3/\langle\alpha_3\alpha_2\alpha_1, \alpha_1\alpha_3\rangle$ \\
          $\mathscr{X}=\emptyset$, so $A$ is not quasi-hereditary by Theorem \ref{main thm}.
\end{description}

\subsection{A Brief Survey on Quasi-Hereditary Nakayama Algebra }

We include a very brief survey in this subsection of quasi-hereditary Nakayama algebras without proof but with some remarks. 

We recall first some notation of Green and Schroll \cite{GS}. Let $Q$ be a finite quiver and let $\mathcal{T}$ be a set of paths in $Q$. A vertex $v$ is said to be {\it properly internal to} $\mathcal{T}$ if some $p\in\mathcal{T}$ has a factorization $p=p_1vp_2$, where $p_1$ and $p_2$ are paths each containing at least one arrow.

Let $A=KQ_n/I$ be  a Nakayama algebra with $I=\langle g_1,\dots, g_k\rangle$.
Let $v_1,\dots,v_m$ be distinct vertices of $Q_n$ and write $v=v_1+\cdots+v_m$. Let $Q_{\hat v}$ be the full subquiver obtained by deleting these vertices, and let $\mathcal{T}_{\hat v}$ consist of those $g_i$ that lie entirely in $Q_{\hat v}$. When no vertex has been deleted, $\mathcal{T}_{\hat 0}=\{g_1,\dots,g_k\}$.

The following theorem collects several existing criteria for quasi-heredity of Nakayama algebras.

\begin{thm}\label{5}
  Let  $A$ be a Nakayama algebra. The following statements are equivalent:

    (1) $A$ is quasi-hereditary.

    (2) There is a simple module of projective dimension $2$, or $A$ is hereditary.

    (3) There is a total ordering of the simple modules $v_1,\dots, v_n$ such that for each $i$, the simple module $v_i$ is not properly internal to $\mathcal{T}_{\widehat{v_1+\dots+v_{i-1}}}$.

    (4) $A$ is $S$-connected.
    
    (5) $\mathscr{X}_A\neq\emptyset$.   
\end{thm}

Conditions (2) and (4) are homological criteria: they require information about the projective dimensions of simple modules. In (4), $S$-connectedness includes finite global dimension, as in \cite{MS}. Condition (3), due to Green and Schroll, is a criterion for monomial algebras formulated in terms of paths in the quiver. In the Nakayama case, condition (5) gives a particularly simple form of this type of criterion: quasi-heredity is determined by the set of simple modules that do not occur in the interiors of the relations.
\vskip 6pt

\section{Applications: Formula for $\q(A)$ and the $q$-Ordering Conjecture}
In this section, we give two applications of the $q$-ordering criterion developed above. One is an iteration formula for $\q(A)$ when $A$ is a Nakayama algebra. The other is a proof of the $q$-ordering conjecture for Nakayama algebras. The section ends with a separate non-monomial example showing that the conjecture does not hold for general finite-dimensional algebras.

\subsection{Formula for $\q(A)$}

We recall the single-relation formula from \cite[Theorem 2]{ZWG} and give a short proof using hoods. 
\begin{thm}( \cite{ZWG}, Theorem 2)
\label{one generator formula}
Let $A=KA_n/I$ be a directed Nakayama algebra with $I=\langle g\rangle$ generated by a path $g$. Then $\q(A)=\frac{2}{\ell_g}n!$.
    \begin{proof}
Since $|\Hod(g)|=\ell_g$, we may assume that $\Hod(g)=\{i,i+1,\dots,i+\ell_g-1\}$. By Theorem \ref{quasi-hereditary order thm}, a total ordering is a $q$-ordering for $A$ if and only if the maximum of $\Hod(g)$ is one of its two endpoints, $i$ and $i+\ell_g-1$. There are two choices for this maximum. Once it is fixed, the remaining $\ell_g-1$ elements of $\Hod(g)$ can be ordered in $(\ell_g-1)!$ ways. The positions of the $\ell_g$ elements of $\Hod(g)$ among the $n$ positions can be chosen in $\binom{n}{\ell_g}$ ways, and the remaining $n-\ell_g$ simple modules can be ordered in $(n-\ell_g)!$ ways. Therefore,
\[
\q(A)=2(\ell_g-1)!\binom{n}{\ell_g}(n-\ell_g)!=\frac{2}{\ell_g}n!.
\]
    \end{proof}
\end{thm}

We next derive an iteration formula for $\q(A)$ by partitioning the $q$-orderings of a Nakayama algebra $A$ according to their maximal simple module. 

Let $I=\langle g_1,\dots,g_k\rangle$ be an admissible ideal of $KQ_n$, with the fixed minimal set of path generators. For $x\in\mathscr{X}_A$, define
\[
I_{\hat{x}}=\langle g_i\mid x\notin\{h(g_i),d(g_i)\},\ 1\leq i\leq k\rangle,
\]
where the empty set generates the zero ideal. Thus $I_{\hat{x}}$ is obtained by deleting every relation having $x$ as an endpoint. If $x\in\mathscr{X}_A^0$, no relation is deleted; if $x\in\mathscr{X}_A^1$, one relation is deleted. If $x\in\mathscr{X}_A^2$, the relations ending and starting at $x$ are deleted, with a relation whose endpoints coincide deleted only once.

Suppose $\preceq$ is a total ordering of $\mathcal{S}$ and $x\in\mathcal{S}$. We use $\max\limits_{\preceq}\mathcal{S}$ to denote the maximum of $\mathcal{S}$ under $\preceq$. Let
\begin{align*}
\mathscr{O}_x(A)=\{\preceq\in\mathscr{O}(A)|x=\max\limits_{\preceq}\mathcal{S}\}
\end{align*}
Write $\q_x(A)=|\mathscr{O}_x(A)|$.

The deletion $I_{\hat{x}}$ can be the zero ideal. This causes no difficulty for a directed quiver, but if $Q_n=\widetilde A_n$, then $K\widetilde A_n$ is infinite-dimensional and $\q(K\widetilde A_n)$ is not defined in our finite-dimensional setting. We therefore isolate this terminal case by setting
\[
\mathscr{Z}_A=\{x\in\mathscr{X}_A\setminus\mathscr{X}_A^0\mid I_{\hat{x}}=0\}
\]
and
\[
\mathscr{Y}_A=(\mathscr{X}_A^1\cup\mathscr{X}_A^2)\setminus\mathscr{Z}_A.
\]

\begin{lem}\label{q_x(A/I)}
Let $A=KQ_n/I$ be a Nakayama algebra with $I=\langle g_1,\dots,g_k\rangle$, and let $x\in\mathcal{S}$. Then
\begin{enumerate}
    \item If $x\notin\mathscr{X}_A$, then $\q_x(A)=0$.
    \item If $x\in\mathscr{X}_A^0$, then $\q_x(A)=\dfrac{1}{n}\q(A)$.
    \item If $x\in\mathscr{Z}_A$, then $\q_x(A)=(n-1)!$.
    \item If $x\in\mathscr{Y}_A$, then
    \[
        \q_x(A)=\frac{1}{n}\q(KQ_n/I_{\hat{x}}).
    \]
\end{enumerate}
    \begin{proof}
For (1), if $x\notin\mathscr{X}_A$, then $x\in\Int(g_i)$ for some $i$. If $x$ were the maximum of $\mathcal{S}$, it would also be the maximum of $\Hod(g_i)$, contrary to Theorem \ref{quasi-hereditary order thm}. Hence $\mathscr{O}_x(A)$ is empty.

For (2), suppose $x\in\mathscr{X}_A^0$. Then $x$ belongs to no hood, so its position in a total ordering does not affect any of the conditions in Theorem \ref{quasi-hereditary order thm}. Fix a relative ordering of the other $n-1$ simple modules that satisfies the hood conditions. Inserting $x$ in any of the $n$ positions then gives a $q$-ordering, and exactly one insertion makes $x$ the maximum. Conversely, every $q$-ordering arises from such a relative ordering. Therefore $|\mathscr{O}_x(A)|=|\mathscr{O}(A)|/n$.

For (3), let $x\in\mathscr{Z}_A$. By the definition of $I_{\hat{x}}$, all generators of $I$ are deleted, and $x$ is an endpoint of every deleted relation. Since $x\in\mathscr{X}_A$, it belongs to none of their interiors. Thus every total ordering having $x$ as its maximum satisfies the criterion of Theorem \ref{quasi-hereditary order thm}. The other $n-1$ simple modules can be ordered arbitrarily, giving $\q_x(A)=(n-1)!$.

For (4), let $x\in\mathscr Y_A$. Every relation deleted in passing to $I_{\hat{x}}$ has $x$ as an endpoint and, since $x\in\mathscr X_A$, does not contain $x$ in its interior. For an ordering whose maximum is $x$, all deleted hood conditions are therefore automatic, while the remaining conditions are unchanged. Consequently,
\[
\mathscr{O}_x(KQ_n/I)=\mathscr{O}_x(KQ_n/I_{\hat{x}}).
\]
After the deletion, $x$ belongs to no remaining hood and hence to the corresponding set $\mathscr X^0$. Since $x\notin\mathscr Z_A$, the quotient remains finite-dimensional. Applying (2) to this algebra gives
\[
\q_x(A)=\q_x(KQ_n/I_{\hat{x}})=\frac{1}{n}\q(KQ_n/I_{\hat{x}}).
\]
    \end{proof}
\end{lem}

We can now prove the following iteration formula for $\q(A)$ via $\mathscr{X}_A$.

\begin{thm}\label{GS decomposition}
Let $A=KQ_n/I$ be a Nakayama algebra. Then
\begin{align}\label{formula q(A)}
\q(A)=\frac{1}{n}\left[|\mathscr{X}_A^0|\q(A)+\sum_{x\in\mathscr{Y}_A}\q(KQ_n/I_{\hat{x}})\right]+|\mathscr{Z}_A|(n-1)!.
\end{align}
\begin{proof}
By partitioning $\mathscr{O}(A)$ according to its maximum and applying Lemma \ref{q_x(A/I)}, we obtain
\begin{align*}
\q(A)
 &=\sum_{x\in\mathcal{S}}\q_x(A)\\
 &=\sum_{x\in\mathscr{X}_A^0}\frac{1}{n}\q(A)
   +\sum_{x\in\mathscr{Y}_A}\frac{1}{n}\q(KQ_n/I_{\hat{x}})
   +\sum_{x\in\mathscr{Z}_A}(n-1)!,
\end{align*}
which is Formula (\ref{formula q(A)}).
\end{proof}
\end{thm}
When $n\neq|\mathscr{X}_A^0|$, Theorem \ref{GS decomposition} can be rewritten as
\[
\q(A)=\frac{\displaystyle\sum_{x\in\mathscr{Y}_A}\q(KQ_n/I_{\hat{x}})+|\mathscr{Z}_A|n!}{n-|\mathscr{X}_A^0|}.
\]
If $I=0$, finite-dimensionality forces $Q_n=A_n$, and $\q(A)=n!$. If $I\neq0$, then $|\mathscr X_A^0|<n$, and each $x\in\mathscr Y_A$ deletes at least one minimal relation while leaving a nonempty set of relations. Thus the formula reduces the number of relations at each recursive step.

The following result is a direct consequence of Theorem \ref{quasi-hereditary order thm}. It also provides the terminal single-relation case for the iteration formula when deleting a relation would leave a cyclic quiver with the zero ideal.

\begin{cor}\label{4.4}

Let $A=KQ_n/\langle g\rangle$ be a Nakayama algebra whose admissible ideal is generated by a single path $g$. Then

\begin{align}\label{onegenerator}
  \q(A)=
    \left\{
\begin{array}{ll}
   \frac{2}{\ell_g}n!,  & \ell_g\leq n; \\
  (n-1)!,  & \ell_g=n+1; \\
  0,  & \ell_g> n+1.
\end{array}
\right.
  \end{align}
\begin{proof}
By Theorem \ref{quasi-hereditary order thm}, a total ordering $\preceq$ of $\mathcal{S}$ is a $q$-ordering if and only if
\[
\max\nolimits_{\preceq}\Hod(g)\notin\Int(g).
\]

Suppose first that $\ell_g\leq n$. Then the $\ell_g$ elements of $\Hod(g)$ are distinct, and the above condition says that the maximum of $\Hod(g)$ must be either $h(g)$ or $d(g)$. There are two choices for this maximum. After it has been fixed, the remaining $\ell_g-1$ elements of $\Hod(g)$ can be ordered in $(\ell_g-1)!$ ways. The positions occupied by the elements of $\Hod(g)$ can be chosen in $\binom{n}{\ell_g}$ ways, and the remaining $n-\ell_g$ simple modules can be ordered in $(n-\ell_g)!$ ways. Hence
\[
\q(A)=2(\ell_g-1)!\binom{n}{\ell_g}(n-\ell_g)!=\frac{2}{\ell_g}n!.
\]

If $\ell_g=n+1$, then $h(g)=d(g)$ and $\Int(g)$ contains every simple module except $h(g)$. Consequently, the condition of Theorem \ref{quasi-hereditary order thm} holds precisely when $h(g)$ is the maximum of $\mathcal{S}$. The remaining $n-1$ simple modules may be ordered arbitrarily, and therefore $\q(A)=(n-1)!$.

Finally, if $\ell_g>n+1$, then every simple module occurs in $\Int(g)$. Thus the maximum of $\Hod(g)$ necessarily belongs to $\Int(g)$, so no total ordering is a $q$-ordering. Hence $\q(A)=0$.
\end{proof}
\end{cor}

Let $A=KQ_n/I$ be a Nakayama algebra with $I=\langle g_1,\dots,g_k\rangle$. Formula (\ref{formula q(A)}) gives a recursive method for computing $\q(A)$, as explained in the following example.

\begin{exm}
    Let $A=KA_5/\langle g_1,g_2\rangle$ with $g_1=\alpha_2\alpha_1$, $g_2=\alpha_3\alpha_2$ be the algebra in Example \ref{A_5}. By Theorem \ref{GS decomposition}, we know
\begin{align*}
    \q(A)&=\frac{1}{5}[\q(A)+\sum\limits_{x\in \mathscr{X}^1}\q(KA_5/I_{\hat{x}})]\\
    &=\frac{1}{5}[\q(A)+\q(KA_5/I_{\hat{1}})+\q(KA_5/I_{\hat{4}})]\\
    &=\frac{1}{4}[\q(KA_5/I_{\hat{1}})+\q(KA_5/I_{\hat{4}})]\\
    &=\frac{1}{4}[\q(KA_5/\langle g_2\rangle)+\q(KA_5/\langle g_1\rangle)]\\
    &=\frac{1}{4}(\frac{2}{3}5!+\frac{2}{3}5!)=\frac{1}{3}5!=40.
\end{align*}
The last line uses Corollary \ref{4.4}.
\end{exm}

\subsection{Sharp bounds and the $q$-ordering conjecture}

For each fixed $n$, we determine the maximum of $\q(A)$ among non-hereditary Nakayama algebras with $n$ simple modules and prove the $q$-ordering conjecture for this class.

The next lemma compares algebras obtained by retaining subsets of a fixed minimal set of relations. Its strict inequality will be used to characterize equality in the upper bound.

\begin{lem}\label{strictly more}
Let $I=\langle g_1,\ldots,g_k\rangle$ be an admissible ideal minimally generated by paths, with $k\geq1$, and put $A=KQ_n/I$. For $1\leq r\leq k$, let $J=\langle g_1,\ldots,g_r\rangle$. Then
\[
\mathscr{O}(A)\subseteq\mathscr{O}(KQ_n/J),
\qquad \q(A)\leq\q(KQ_n/J).
\]
Moreover, if $n\geq3$, $\ell_{g_1}=3$, and $k\geq2$, then
\[
\q(A)<\q(KQ_n/\langle g_1\rangle)=\frac{2}{3}n!.
\]
\begin{proof}
The retained paths form a minimal generating set of $J$. The quotient $KQ_n/J$ is finite-dimensional: this is automatic in the directed case, and in the cyclic case every sufficiently long path contains a retained relation. Since a $q$-ordering for $A$ satisfies the hood condition for every $g_i$, it satisfies those for $g_1,\ldots,g_r$. The inclusion and the weak inequality follow from Theorem \ref{quasi-hereditary order thm}.

For the strict assertion, let $v$ be the unique interior vertex of $g_1$. The hood of $g_1$ has three distinct vertices because $n\geq3$. We claim that $\Int(g_2)$ contains a vertex $y\neq v$. If $\ell_{g_2}=3$, equality of the interior vertices would force $g_2=g_1$, since a path with two arrows in $Q_n$ is determined by its middle vertex. This contradicts minimality. If $\ell_{g_2}>3$, its interior contains two consecutive, distinct vertices, at least one of which differs from $v$.

Choose a total ordering in which $y$ is the maximum of $\mathcal S$ and the maximum of $\Hod(g_1)$ is an endpoint. Such an ordering exists: if $y$ belongs to $\Hod(g_1)$, then it is already an endpoint; otherwise, place $h(g_1)$ above the other vertices of that hood. By Theorem \ref{quasi-hereditary order thm}, this ordering belongs to $\mathscr{O}(KQ_n/\langle g_1\rangle)$, but it does not belong to $\mathscr{O}(A)$ because $y\in\Int(g_2)$. Thus the inclusion is strict. Corollary \ref{4.4} gives $\q(KQ_n/\langle g_1\rangle)=\frac{2}{3}n!$.
\end{proof}
\end{lem}

\begin{thm}\label{2/3}
Let $A=KQ_n/I$ be a non-hereditary Nakayama algebra with $I$ an admissible ideal of $KQ_n$. Then
  \begin{align}\q(A)\le \frac{2}{3}n!
    \end{align}
    Equality holds if and only if $n\geq3$ and $I$ is generated by a single path of length $3$ (that is, a path with two arrows). When $n=1$, $\q(A)=0$; when $n=2$, one has the sharper bound $\q(A)\leq1$. In particular, the $q$-ordering conjecture is true for Nakayama algebras.
    \begin{proof}
Since $A$ is non-hereditary, $I\neq0$. Write $I=\langle g_1,\ldots,g_k\rangle$ minimally, with $k\geq1$. By Lemma \ref{strictly more},
\[
\q(A)\leq\q(KQ_n/\langle g_i\rangle)
\]
for each $i$, after relabelling the generators if necessary.

If $n=1$, then $\ell_{g_i}\geq3>n+1$, so Corollary \ref{4.4} gives $\q(A)=0<\frac{2}{3}$.

If $n=2$, the quiver must be cyclic, because the directed quiver $A_2$ has no paths with two arrows and admits no nonzero admissible ideal. Since $\ell_{g_i}\geq3=n+1$, Corollary \ref{4.4} gives $\q(KQ_n/\langle g_i\rangle)\leq1$. Hence $\q(A)\leq1<\frac{2}{3}\,2!$. The value $1$ is attained by $K\widetilde A_2/\langle\alpha_2\alpha_1\rangle$.

Now let $n\geq3$. For a single relation $g$, Corollary \ref{4.4} gives
\[
\frac{\q(KQ_n/\langle g\rangle)}{n!}=
\begin{cases}
2/\ell_g,&3\leq\ell_g\leq n,\\
1/n,&\ell_g=n+1,\\
0,&\ell_g>n+1.
\end{cases}
\]
Each value is at most $2/3$, with equality precisely when $\ell_g=3$. This proves the upper bound for $A$. If $\q(A)=\frac{2}{3}n!$, comparison with each single-relation quotient forces every $g_i$ to have length $3$. If $k\geq2$, the strict part of Lemma \ref{strictly more} gives a contradiction. Thus $k=1$. Conversely, when $n\geq3$ and $I$ is generated by a single relation of length $3$, Corollary \ref{4.4} gives equality.
    \end{proof}
\end{thm}

Theorem \ref{2/3} shows that the $q$-ordering conjecture holds for all Nakayama algebras. The following non-monomial example shows that the conjecture fails for general finite-dimensional algebras.

\begin{exm}
Let $A=KQ/I$, where $I=\langle\alpha_2\alpha_1-\alpha_4\alpha_3\rangle$.
\[
\xymatrix{
& & &  *+{2} \ar[dr]^{\alpha_2}\\
&Q: & *+{1} \ar[ur]^{\alpha_1} \ar[dr]_{\alpha_3} & &*+{4} \\
& & & *+{3} \ar[ur]_{\alpha_4}  \\
}
\]
Only the four orderings
\[
1 \preceq 4 \preceq 2 \preceq 3,\quad
1 \preceq 4 \preceq 3 \preceq 2,\quad
4 \preceq 1 \preceq 2 \preceq 3,\quad
4 \preceq 1 \preceq 3 \preceq 2
\]
are not $q$-orderings. Thus, $\q(A)=20>\frac{2}{3}\cdot4!$.
\end{exm}

\section*{Declarations}
\subsection*{Conflict of interest} The authors state that there is no conflict of interest.

\subsection*{Funding}
This research is supported by  National Natural Science Foundation of China under Grant No. 123B1026.

\subsection*{Ethical approval}
This article does not contain any studies with human participants or animals performed by any of the authors.

\subsection*{Informed consent}
Not applicable.

\subsection*{Author Contributions}
Both authors contributed equally to this work.

\subsection*{Data Availability Statement}
This article does not involve any data, as it is a purely theoretical study.

\bibliographystyle{plain}

\bibliography{main}

@article {GS,
    AUTHOR = {Green, Edward L. and Schroll, Sibylle},
     TITLE = {On quasi-hereditary algebras},
   JOURNAL = {Bull. Sci. Math.},
  FJOURNAL = {Bulletin des Sciences Math\'{e}matiques},
    VOLUME = {157},
      YEAR = {2019},
     PAGES = {102797, 14},
      ISSN = {0007-4497,1952-4773},
   MRCLASS = {16G20 (16D10 16G10)},
  MRNUMBER = {4017924},
MRREVIEWER = {Karin\ M.\ Jacobsen},
       DOI = {10.1016/j.bulsci.2019.102797},
       URL = {https://doi.org/10.1016/j.bulsci.2019.102797},
}

@article {UY,
    AUTHOR = {Uematsu, Morio and Yamagata, Kunio},
     TITLE = {On serial quasi-hereditary rings},
   JOURNAL = {Hokkaido Math. J.},
  FJOURNAL = {Hokkaido Mathematical Journal},
    VOLUME = {19},
      YEAR = {1990},
    NUMBER = {1},
     PAGES = {165--174},
      ISSN = {0385-4035},
   MRCLASS = {16E60 (16E10)},
  MRNUMBER = {1039471},
MRREVIEWER = {Peter\ Dr\"{a}xler},
       DOI = {10.14492/hokmj/1381517167},
       URL = {https://doi.org/10.14492/hokmj/1381517167},
}

@article {MS,
    AUTHOR = {Marczinzik, Ren\'{e} and Sen, Emre},
     TITLE = {A new characterization of quasi-hereditary {N}akayama algebras
              and applications},
   JOURNAL = {Comm. Algebra},
  FJOURNAL = {Communications in Algebra},
    VOLUME = {50},
      YEAR = {2022},
    NUMBER = {10},
     PAGES = {4481--4493},
      ISSN = {0092-7872,1532-4125},
   MRCLASS = {16G10 (11B39 16E10)},
  MRNUMBER = {4447471},
MRREVIEWER = {Liping\ Li},
       DOI = {10.1080/00927872.2022.2063301},
       URL = {https://doi.org/10.1080/00927872.2022.2063301},
}

@article {CPS,
    AUTHOR = {Cline, E. and Parshall, B. and Scott, L.},
     TITLE = {Finite-dimensional algebras and highest weight categories},
   JOURNAL = {J. Reine Angew. Math.},
  FJOURNAL = {Journal f\"{u}r die Reine und Angewandte Mathematik. [Crelle's
              Journal]},
    VOLUME = {391},
      YEAR = {1988},
     PAGES = {85--99},
      ISSN = {0075-4102,1435-5345},
   MRCLASS = {18E30 (16A46 17B10 20G05 32C38)},
  MRNUMBER = {961165},
MRREVIEWER = {Dieter\ Happel},
}

@article {ZL,
    AUTHOR = {Zhang, Yue Hui and Li, Yu},
     TITLE = {Counting Quasi-Hereditary Orderings of Finite Dimensional Algebras},
   JOURNAL = {J. MATH. TECH.},
  FJOURNAL = {JOURNAL OF MATHEMATICS FOR TECHNOLOGY},
    VOLUME = {16},
      YEAR = {2000},
    NUMBER = {3},
     PAGES = {9--11},
      ISSN = {1672-1454},
   MRCLASS = {16G10},
  MRNUMBER = {=},
MRREVIEWER = {=},
}

@article {ZWG,
    AUTHOR = {Zhang, Yue Hui and Wu, Liu San and Gao, Chun Yan},
     TITLE = {Quasi-hereditary orderings of {$A_n$}-type algebras with two
              generators},
   JOURNAL = {J. Math. Res. Exposition},
  FJOURNAL = {Journal of Mathematical Research and Exposition},
    VOLUME = {28},
      YEAR = {2008},
    NUMBER = {4},
     PAGES = {975--980},
      ISSN = {1000-341X},
   MRCLASS = {16G10},
  MRNUMBER = {2465212},
MRREVIEWER = {Rita\ Zuazua},
}

@article {C,
    AUTHOR = {Coulembier, Kevin},
     TITLE = {The classification of blocks in {BGG} category {$\mathcal{O}$}},
   JOURNAL = {Math. Z.},
  FJOURNAL = {Mathematische Zeitschrift},
    VOLUME = {295},
      YEAR = {2020},
    NUMBER = {1-2},
     PAGES = {821--837},
      ISSN = {0025-5874,1432-1823},
   MRCLASS = {17B10 (17B20 18E10 20F55)},
  MRNUMBER = {4100034},
MRREVIEWER = {Luis\ Enrique\ Ramirez},
       DOI = {10.1007/s00209-019-02376-9},
       URL = {https://doi.org/10.1007/s00209-019-02376-9},
}

@incollection {DR1992,
    AUTHOR = {Dlab, Vlastimil and Ringel, Claus Michael},
     TITLE = {The module theoretical approach to quasi-hereditary algebras},
 BOOKTITLE = {Representations of algebras and related topics ({K}yoto,
              1990)},
    SERIES = {London Math. Soc. Lecture Note Ser.},
    VOLUME = {168},
     PAGES = {200--224},
 PUBLISHER = {Cambridge Univ. Press, Cambridge},
      YEAR = {1992},
      ISBN = {0-521-42411-9},
   MRCLASS = {16E60 (16D90 16G20)},
  MRNUMBER = {1211481},
MRREVIEWER = {Otto\ Kerner},
}

@article {FKR,
    AUTHOR = {Flores, Manuel and Kimura, Yuta and Rognerud, Baptiste},
     TITLE = {Combinatorics of quasi-hereditary structures},
   JOURNAL = {J. Combin. Theory Ser. A},
  FJOURNAL = {Journal of Combinatorial Theory. Series A},
    VOLUME = {187},
      YEAR = {2022},
     PAGES = {Paper No. 105559, 54},
      ISSN = {0097-3165,1096-0899},
   MRCLASS = {16E60 (05E16)},
  MRNUMBER = {4343742},
MRREVIEWER = {Sudalaimuthu\ Santhakumar},
       DOI = {10.1016/j.jcta.2021.105559},
       URL = {https://doi.org/10.1016/j.jcta.2021.105559},
}

@article {G,
    AUTHOR = {Goto, Yuichiro},
     TITLE = {Criterion for quasi-heredity},
   JOURNAL = {Algebr. Represent. Theory},
  FJOURNAL = {Algebras and Representation Theory},
    VOLUME = {27},
      YEAR = {2024},
    NUMBER = {2},
     PAGES = {1395--1403},
      ISSN = {1386-923X,1572-9079},
   MRCLASS = {16G20 (06A05 20B30)},
  MRNUMBER = {4741527},
       DOI = {10.1007/s10468-024-10263-z},
       URL = {https://doi.org/10.1007/s10468-024-10263-z},
}

@article {zp1,
    AUTHOR = {Zhang, Pu},
     TITLE = {Quasi-hereditary algebras and tilting modules},
   JOURNAL = {Comm. Algebra},
  FJOURNAL = {Communications in Algebra},
    VOLUME = {24},
      YEAR = {1996},
    NUMBER = {12},
     PAGES = {3707--3717},
      ISSN = {0092-7872,1532-4125},
   MRCLASS = {16G10},
  MRNUMBER = {1408501},
MRREVIEWER = {Steffen\ K\"onig},
       DOI = {10.1080/00927879608825781},
       URL = {https://doi.org/10.1080/00927879608825781},
}

@article {zp2,
    AUTHOR = {Zhang, Pu},
     TITLE = {A criterion for quasi-heredity and the characteristic module},
   JOURNAL = {Manuscripta Math.},
  FJOURNAL = {Manuscripta Mathematica},
    VOLUME = {93},
      YEAR = {1997},
    NUMBER = {2},
     PAGES = {129--135},
      ISSN = {0025-2611,1432-1785},
   MRCLASS = {16G10},
  MRNUMBER = {1464360},
MRREVIEWER = {Lidia\ Angeleri H\"ugel},
       DOI = {10.1007/BF02677461},
       URL = {https://doi.org/10.1007/BF02677461},
}

@article {ringer91,
    AUTHOR = {Ringel, Claus Michael},
     TITLE = {The category of modules with good filtrations over a
              quasi-hereditary algebra has almost split sequences},
   JOURNAL = {Math. Z.},
  FJOURNAL = {Mathematische Zeitschrift},
    VOLUME = {208},
      YEAR = {1991},
    NUMBER = {2},
     PAGES = {209--223},
      ISSN = {0025-5874,1432-1823},
   MRCLASS = {16G10 (16D90)},
  MRNUMBER = {1128706},
MRREVIEWER = {J.\ Antonio\ de la Pe\~na},
       DOI = {10.1007/BF02571521},
       URL = {https://doi.org/10.1007/BF02571521},
}

@article {DR89,
    AUTHOR = {Dlab, Vlastimil and Ringel, Claus Michael},
     TITLE = {Quasi-hereditary algebras},
   JOURNAL = {Illinois J. Math.},
  FJOURNAL = {Illinois Journal of Mathematics},
    VOLUME = {33},
      YEAR = {1989},
    NUMBER = {2},
     PAGES = {280--291},
      ISSN = {0019-2082,1945-6581},
   MRCLASS = {16A35 (16A48 16A60)},
  MRNUMBER = {987824},
MRREVIEWER = {Alfred\ G.\ Wiedemann},
       URL = {http://projecteuclid.org/euclid.ijm/1255988725},
}

@article {K89,
    AUTHOR = {K\"onig, Steffen and Wiedemann, Alfred},
     TITLE = {Global dimension two orders are quasi-hereditary},
   JOURNAL = {Manuscripta Math.},
  FJOURNAL = {Manuscripta Mathematica},
    VOLUME = {66},
      YEAR = {1989},
    NUMBER = {1},
     PAGES = {17--23},
      ISSN = {0025-2611,1432-1785},
   MRCLASS = {16A18},
  MRNUMBER = {1023548},
MRREVIEWER = {Wolfgang\ Rump},
       DOI = {10.1007/BF02568479},
       URL = {https://doi.org/10.1007/BF02568479},
}

@article {MRS,
    AUTHOR = {Marczinzik, Ren\'e{} and Rubey, Martin and Stump, Christian},
     TITLE = {A combinatorial classification of 2-regular simple modules for
              {N}akayama algebras},
   JOURNAL = {J. Pure Appl. Algebra},
  FJOURNAL = {Journal of Pure and Applied Algebra},
    VOLUME = {225},
      YEAR = {2021},
    NUMBER = {3},
     PAGES = {Paper No. 106520, 34},
      ISSN = {0022-4049,1873-1376},
   MRCLASS = {16G10 (05A15 05A19 16E10)},
  MRNUMBER = {4142195},
MRREVIEWER = {Chao\ Zhang},
       DOI = {10.1016/j.jpaa.2020.106520},
       URL = {https://doi.org/10.1016/j.jpaa.2020.106520},
}

@article{RINGEL202195,
title = {The finitistic dimension of a {Nakayama} algebra},
journal = {Journal of Algebra},
volume = {576},
pages = {95-145},
year = {2021},
issn = {0021-8693},
doi = {https://doi.org/10.1016/j.jalgebra.2021.01.040},
url = {https://www.sciencedirect.com/science/article/pii/S0021869321000909},
author = {Claus Michael Ringel}
}

@article{Sen03062021,
author = {Emre Sen},
title = {The $\phi$-dimension of cyclic {Nakayama} algebras},
journal = {Communications in Algebra},
volume = {49},
number = {6},
pages = {2278--2299},
year = {2021},
publisher = {Taylor \& Francis},
doi = {10.1080/00927872.2020.1870995},


URL = { 
    
        https://doi.org/10.1080/00927872.2020.1870995
    
    

},
eprint = { 
    
        https://doi.org/10.1080/00927872.2020.1870995
    
    

}

}

@article{RINGEL2013241,
title = {The {Gorenstein} projective modules for the {Nakayama} algebras. {I}},
journal = {Journal of Algebra},
volume = {385},
pages = {241-261},
year = {2013},
issn = {0021-8693},
doi = {https://doi.org/10.1016/j.jalgebra.2013.03.014},
url = {https://www.sciencedirect.com/science/article/pii/S0021869313001609},
author = {Claus Michael Ringel}
}
\end{document}